INSTITUT NATIONAL DE RECHERCHE EN INFORMATIQUE ET EN AUTOMATIQUE

# *The Genetic Code Degeneration I: Rules Governing the Code Degeneration and the Spatial Organization of the Codon Informative Properties.*

Melina Rapacioli, Edmundo Rofman, Vladimir Flores

Rapport de recherche

# The Genetic Code Degeneration I: Rules Governing the Code Degeneration and the Spatial Organization of the Codon Informative Properties.

Melina Rapacioli [1], Edmundo Rofman [2], Vladimir Flores [1]

[1] Interdisciplinary Group in Theoretical Biology, Department of Biostructural Sciences, Favaloro University. Solís 453 (1078) Buenos Aires, Argentina.
vflores@favaloro.edu.ar

[2] Inria-Rocquencourt, France and Instituto Argentino de Matematicas (IAM-CONICET), Argentine.
Edmundo.Rofman@inria.fr

**Abstract:** the present work is devoted to describe a set of rules explaining the discriminating versus non-discriminating behavior of the di-basic stages and to characterize the role of each base in determining such a behavior. Bases are analyze as dual entities characterized by its chemical type and the number of H bonds involved in the codon-anticodon interaction. A codon is characterized as an asymmetric informative entity whose global informative capacity results from the spatially organized combinatory of the 6 properties assigned by the 3 bases.

This work was supported by grants from CONICET, Argentina

# La dégénération du code génétique I: les règles qui gouvernent la dégénération du code et l´organisation spatiale des propriétés informatives du codon.

**Résumé:** Ce travail a le but de décrire un ensemble de règles qui expliquent le caractère discriminant versus non-discriminant des états di-basiques et de caractériser le rôle de chaque base dans la détermination du caractère déjà mentionné. Les bases sont analysées comme des entités duales, caractérisées par sa nature chimique et par le nombre de ponts de Hydrogène involucrés dans l´interaction codon-anticodon. Le codon est caractérisé comme une entité informative asymétrique, sa capacité informative globale résulte de la combinaison organisée dans l'espace des six propriétés qui lui sont assignées par les trois bases qui le composent.

# 1 Introduction

The genetic information required for protein synthesis is encoded in specific DNA segments, the so called structural genes. The informative value resides on the specific sequence of four nucleotides possessing each of them a distinct nitrogenated base. The specific sequence of bases determines the specific amino acids (AAs) sequence of each single protein. The genetic code refers to the set of rules that specifies how the DNA base sequences of the structural genes are translated into the amino acid sequences of the proteins. These rules are the net result of specific matching: **(a)** between codons (CD) (a set of three adjacent bases) of the mRNA and anti-codons (Anti-CD) of the t-RNA; **(b)** between chemical groups of the t-RNA and specific enzymes (aminoacyl-tRNA synthetases) and **(c)** between these enzymes and the AAs of a nascent protein.

The genetic code is said to be degenerated. Excluding the 3 Stop CDs, the notion of genetic code degeneration refers to the fact that 61 CDs are used to specify 20 amino acids. More specifically the terms refers to the facts that (a) the Anti-CD of some tRNA matches more than one CD and that (b) some AAs can be specified by more than one tRNA type having different Anti-CDs [1, 2].

CDs that specify the same AA are called synonyms. As a rule they share the bases located at the 1$^{st}$ and 2$^{nd}$ position. That is why several authors prefer a CD description as the association between a di-nucletide (those located in 1$^{st}$ and 2$^{nd}$ position) and the nucleotide of the 3$^{rd}$ position [3, 4, 5, 6, 7, 8]. We will use the term di-base instead of di-nucleotide.

The set of 4 CDs corresponding to a particular di-base behave either as a quadruplet [4pts] (one set of 4 synonyms) or as diplets [2pts] (two sets composed of 2 synonyms each) or display an even more complex behavior that will be later described. In the present analysis di-bases will be characterized as discriminating (**D**) or non-discriminating (**non-D**). A di-base will be defined as discriminating when its corresponding 4 CDs diverge into two 2pts depending on the type of base (pyrimidine or purine) of the 3$^{rd}$ position. A di-nucleotide will be defined as non-discriminating when its corresponding 4 CDs conform one 4pt independently of the type of base of the 3$^{rd}$ position.

This work possesses three different approaches related to the genetic code degeneration. The first part, presented in this paper, aims at defining a set of rules that formally describe the discriminating and non-discriminating character of the different kind of di-bases. This paper also tries to characterize the role the different bases play in assigning such characters to di-bases. Bases will be analyzed as dual physicochemical entities possessing two essential characteristics: (a) the chemical or molecular type and (b) the number of Hydrogen bonds involved in the CD-AntiCD interaction. A second part of this work will introduce a novel 2-dimensional space of CD representation based on their dual character. The third will present statistical analyses of CD frequency in mRNA of related families of human proteins and a comparison with non-coding DNA sequences.

# 2 The Rules that Govern the Genetic Code Degeneration

This section introduces enunciations that describe the “logic” implicit in the genetic code and precisely describe the genetic code degeneration. These rules can be deduced from established data about (a) relationship between CDs and AAs and (b) the occurrence of redundant CDs [1, 2, 9].

## 2.1 Definitions and nomenclature

Some terms will be defined and the relationship with current concepts will be clarified.

- Codon [CD] is the informative unit of the mRNA. It is a sequence of 3 adjacent bases [B]. There are four different Bs [uracil (U), cytosine (C), adenine (A) and guanine (G)].

- Each CD is denoted by the sequence of three Bs (BBB). The position of a B within the CD is denoted as: B1 or B- - for the 1st position; B2 or –B– for the 2nd position and B3 or ––B for the 3rd position.

- A di-base is denoted as B1B2 or BB–. There exist 16 BB–.

- For each particular BB– there are 4 CDs given that ––B could be U, C, A or G.

- A BB– is designated as **non-Discriminating** (**non-D**) when the four CDs specify a unique AA. Then, the four CDs of a **non-D** BB– are synonyms forming a 4pt.

- A BB– is designated as **Discriminating** (**D**) when two CDs specify a particular AA and the others two specify another AA. Then, the four CDs conforms 2 groups composed of 2 synonyms each; i.e. they are grouped into two 2pts.

- Each B is characterized by two essential properties: molecular type [mT] and number of H bonds [nHb].

- Two Bs, U and C, belongs to the Pyrimidine type and are denoted **Y**. The others two Bs, A and G, belongs to the Purine type and are denoted **R**.

- Bs belonging to both molecular type (Y and R) may possess 2 or 3 H bonds. This property is denoted as **2** or **3**.

- Two Bs, U and A share the property (nHb) **2**. The others two, C and G, share the property (nHb) **3**.

- Each B may be denoted indicating its specific pair of properties: **Y** or **R** / **2** or **3**.

- Bs of pyrimidine type are denoted as follows: U is $\frac{Y}{2}$; C is $\frac{Y}{3}$.

- Bs of purine type are denoted as follows: A is $\frac{R}{2}$; G is $\frac{R}{3}$.

- A B is defined as **coherent** when its properties (**Y** or **R** and **2** or **3**) coherently contribute to the discriminating or non-discriminating behavior of the BB– they compose.

- A B is defined as **non-coherent** when its properties (**Y** or **R** and **2** or **3**) do not contribute coherently to the discriminating or non-discriminating behavior of the BB– they compose. The notions of **coherence** and **non-coherence** are defined in following paragraphs.

- A CD can be denoted by specifying the 6 properties of its corresponding three Bs

Example: The CD UGC is denoted as $\frac{Y}{2}\frac{R}{3}\frac{Y}{3}$

## 2.2 Rules of the genetic code degeneration

It is known that there are set of synonyms of different sizes (Figure 1): singlets ([1pt]: unitary set of CDs), diplets ([2pt]: a set of 2 synonyms), triplets ([3pt]: set of 3 synonyms), quadruplets ([4pt]: set of 4 synonyms) and sextuplets ([6pt]: set of 6 synonyms). A 3pt can be considered as the association of a 2pt + a 1pt and, a 6pt can be considered as association of a 4pt + a 2pt.
In the simplest diagram (Figure 2), only the discriminating or non-discriminating character of each of the 16 BB– is represented. In order to make the description easier, as a first approach, the 4 CDs corresponding to each BB– will be described as if they were just 4pts or 2pts. In fol-lowing paragraphs the details will be considered step by step, from the simplest picture to the most complex one, and rules of increasing complexity will then be introduced.

| First Position | Second Position: U | Second Position: C | Second Position: A | Second Position: G | Third Position |
|---|---|---|---|---|---|
| **U** | UUU Phe<br>UUC Phe<br>UUA Leu<br>UUG Leu | UCU Ser<br>UCC Ser<br>UCA Ser<br>UCG Ser | UAU Tyr<br>UAC Tyr<br>UAA Stop<br>UAG Stop | UGU Cys<br>UGC Cys<br>UGA Stop<br>UGG Trp | **U**<br>**C**<br>**A**<br>**G** |
| **C** | CUU Leu<br>CUC Leu<br>CUA Leu<br>CUG Leu | CCU Pro<br>CCC Pro<br>CCA Pro<br>CCG Pro | CAU His<br>CAC His<br>CAA Gln<br>CAG Gln | CGU Arg<br>CGC Arg<br>CGA Arg<br>CGG Arg | **U**<br>**C**<br>**A**<br>**G** |
| **A** | AUU Ile<br>AUC Ile<br>AUA Ile<br>AUG Met | ACU Thr<br>ACC Thr<br>ACA Thr<br>ACG Thr | AAU Asn<br>AAC Asn<br>AAA Lys<br>AAG Lys | AGU Ser<br>AGC Ser<br>AGA Arg<br>AGG Arg | **U**<br>**C**<br>**A**<br>**G** |
| **G** | GUU Val<br>GUC Val<br>GUA Val<br>GUG Val | GCU Ala<br>GCC Ala<br>GCA Ala<br>GCG Ala | GAU Asp<br>GAC Asp<br>GAA Glu<br>GAG Glu | GGU Gly<br>GGC Gly<br>GGA Gly<br>GGG Gly | **U**<br>**C**<br>**A**<br>**G** |

**Figure 1. The genetic code. The table illustrates the sets of CDs that specify each amino acid.**

### 2.2.1 Bases A and C

With regards to the bases A and C, it can be seen in Figure 2 that:

> **Rule 1.** Any BB– which contains C is non-discriminating except when B2 is A and any BB– which contains A is discriminating except when B2 is C.

This rule reveals three essential facts:

(1) C assigns a non-D behavior to the BB–. C possesses a non-discriminating character.

(2) A assigns a D behavior to the BB–. A possesses a discriminating character.
BB– composed of combinations of A and C, such as C**A**– and A**C**–, consist of combinations of opposing characteristics. This opposition of characters is solved in a way indicating that B1 and B2 (the 1st and the 2nd positions) display different informative roles.

(3) B2, the 2nd position, has a relevant informative value. In considering those BB– that are combinations of C and A, the discriminating or non-discriminating character is determined by **B2**, i.e. it depends on the character of the B located in the 2nd position.

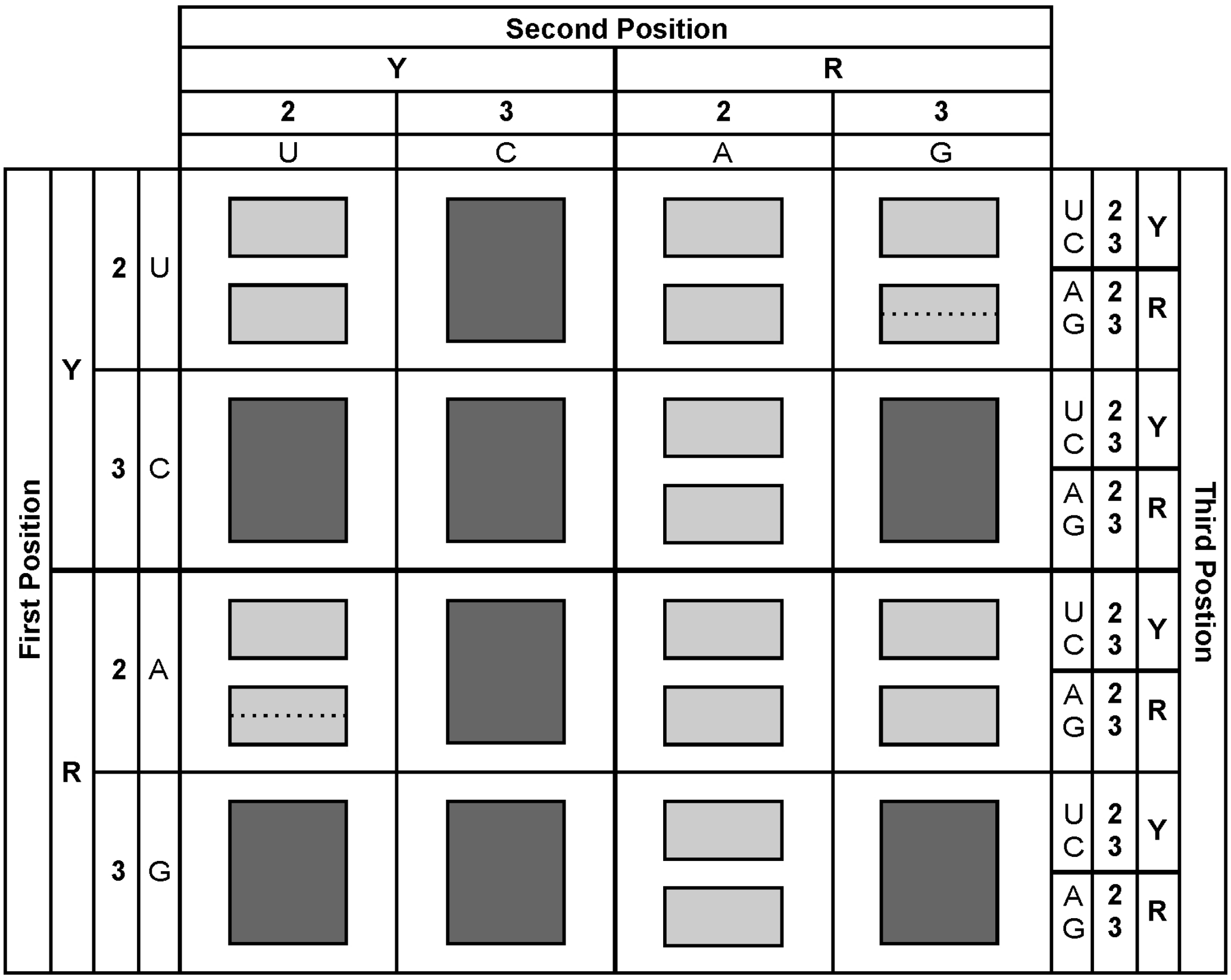

**Figure 2. Simplified picture of the genetic code degeneration. The 4 CDs of each BB- are represented as if they were only 4pts or 2pts. Light gray: 2pts. Dark gray: 4pts. Dotted lines: indicate the divergence of 2pts into 2 1pts. The distribution of the properties Y, R, 2 and 3 is indicated.**

### 2.2.2 Bases U and G

There is no equivalent rule applicable to Bs U and G. It will be explained later that this is due to the fact that these Bs are different from A and C with regards the character defined as coherence.

With regards to BB– composed of combinations of U and G: (a) CDs containing the combination G**U**– corresponds to a 4pt and (b) CDs containing the combination U**G**– diverge in two 2pts.

Then, the structure –**U**– assigns a non-D behavior to a BB– and the structure –**G**– assigns a D character. As a consequence, U possesses a non-Discriminating character and G possesses a discriminating character.

### 2.2.3 The role of the molecular type (Y, R)

From the previous paragraphs can be enunciated another general rule:

**Rule 2. Y** is a **non-discriminating** property and **R** is a **discriminating** property

### 2.2.4 The role of the property nHb (2, 3)

C and G possess the property **3** and U and A are characterized by the property **2**. Analyzing the distribution of 4pts and 2pts, in Figure 2, in connection with the properties **2** and **3**, it can be postulated another general rule:

**Rule 3.** Any BB– whose bases in 1st and 2nd position have the property **3** is non-discriminating and any BB– whose bases in 1st and 2nd position have the property **2** is discriminating.

A more general enunciation of this rule is:

**Rule 4. 3** is a **non-discriminating** property and **2** is a **discriminating** property

All this statements are included in the following **rule**:

Any $-\ \overset{Y}{-}\ -$ belongs to a 4pt except in the case $\underset{2}{-}\ \overset{Y}{\underset{2}{-}}\ -$ and any $-\ \overset{R}{-}\ -$ belongs to a 2pt except in the case $\underset{3}{-}\ \overset{R}{\underset{3}{-}}\ -$

## 2.3 Towards more general rules

In next paragraphs we will present two set of rules implying that the CD structure possesses a dual informative function, on one hand it determines (a) the discriminating or non-discriminating character and, on the other, it determines (b) a specific AA.

### 2.3.1 Rules of determination of the Discriminating or non-Discriminating character of the CD (4pts vs 2pts)

The behavior described by **Rule 1** is due to the fact that: **(a)** base C is characterized by the coincidence of two non-discriminating properties (Y and 3) and **(b)** base A is defined by the coincidence of two discriminating properties (R and 2). We define this condition as **Coherence**. C and A are coherent bases.
Bs U and G do not fulfill this condition. The base U is characterized by one discriminating property associated to a non-discriminating one (**2** and **Y** respectively) and G is defined by the association of one non-discriminating and one discriminating property (**3** and **R** respectively). We define this condition as **non-Coherence**. U and G are non-coherent bases.

The relevance of the spatial position is reveal by the fact that all of the rules previously enunciated can be integrated in a unique enunciation:

**First law**: Three properties located at defined positions are involved in the determination of the discriminating or non-discriminating character of a CD.

$$___ \; \frac{mT_{Y-R}}{nHb_{3-2}nHb_{3-2}} \; ___$$

The other three properties are absolutely irrelevant in this respect. They do not possess informative value with regards to this characteristic of the CD.

**Second Law**: A CD is discriminating, respectively non-discriminating, whenever two of these properties are discriminating, respectively non-discriminating.

A derivation of the $2^{nd}$ law dictates that

A coherent B2 by itself determines the D or non-D character of the CD.

From this, two additional enunciations derive:

Any $- \frac{Y}{3} -$ is non-discriminating, that is belongs to a 4pt, and any $- \frac{R}{2} -$ is discriminating, that is belongs to a 2pt.

The remaining four properties are informatively irrelevant

Another derivation of the first law establishes that

In any CD with a non-coherent B2 the discriminating or non-Discriminating character depends, also on the property nHb of B1.

### 2.3.2 Rules about the AA specifying character of the CD

There exist BB– whose 4 CDs correspond to a 4pt. They specify a single AA.

There exist BB– whose 4 CDs diverge into two 2pts that specify two different AAs.

There exist BB– whose 4 CDs diverge into one 2pt and two 1pts. In some cases the 2pt + one 1pt conform a 3pt that specify a single AA and the 1pt specifies another AA. In other cases the 1pts correspond to Start or Stop CDs.

The following enunciations aim at demonstrating that for each of the above mentioned conditions there are different rules explaining how a CD specify a particular AA. These rules describe that: (a) the “quantity” of information (number of informative elements, i.e. base properties) stored by CDs exceeds that required for the specification of an AA by 4pts; (b) the “quantity” of information (number of informative elements) of the CD required for the specification of an AA by 2pts is a little higher; (c) only when an AA is specify by a 1pt the complete set of informative elements of the CD is required.

**Rule of specification of AAs by 4pts**

**First rule:** Four of the six properties of a CD are required to specify an AA by 4pts.

1. In non-discriminating BB– the AA specification depends on B2+B1

1.1. In non-discriminating BB– with coherent B2: {B2} specifies a 4pt and {B1B2} specifies the AA.

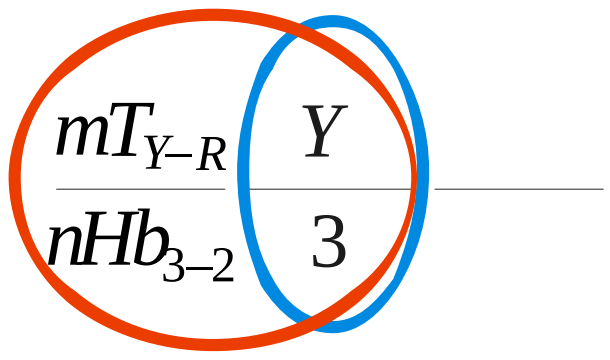

1.2. In non-discriminating BB– with non-coherent B2: {3 in B1+B2} specifies a 4pt and {B1B2} specifies the AA

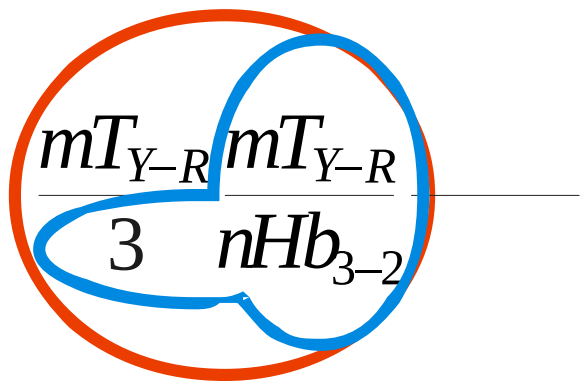

**Rule of AAs specification by 2pts**

**Second rule:** Five of the six properties of a CD are required for the AA specification by 2pts

2. In discriminating BB– the ability to specify an AA depends, additionally, on the property **Y** vs **R** of B3.

2.1. In discriminating BB– with coherent B2: {B2} determines the divergence in two 2pts, {B1B2} specify which pair of 2pts and {B1B2 + mT of B3} specifies a particular AA

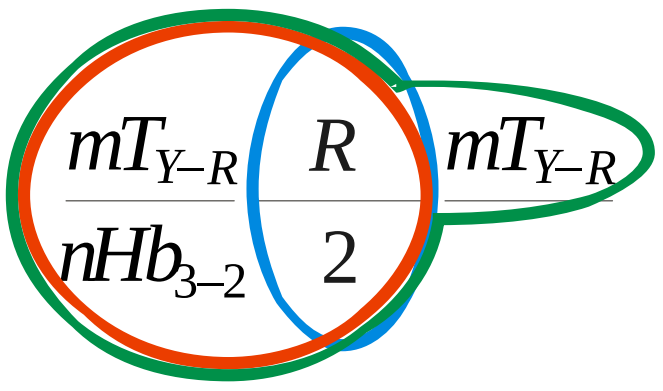

2.2. In discriminating BB– with non-coherent B2: {2 in B1+B2} determine the divergence in two 2pts, {B1B2} specify which pair of 2pts and {B1B2+ mT of B3} specify a particular AA

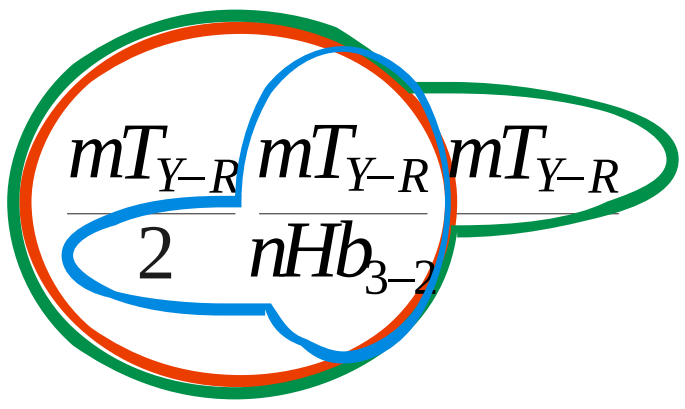

**Rules of AAs specification by 1pts**

**Third rule:** The six properties of a CD are needed for the AA specification by 1pts

The specification of an AA by 1pts can be considered special cases in which a 2pt diverges into two 1pts. This condition requires using the complete set or properties of the CD and takes place only when a discriminating BB– possesses a non-coherent B2 and R in B3. In these cases the AA specification depends also on the properties 2 vs 3 of B3.

In discriminating BB– with non-coherent B2 and R in B3, the BB– can diverge originating 1pts: {2 in B1+ B2} determine the discriminating character, {B1B2} identify a particular BB–, {B1B2+ R in B3} allows the divergence to 1pts, {B1B2+ R in B3+ 2 vs 3 in B3} specify a particular AA.

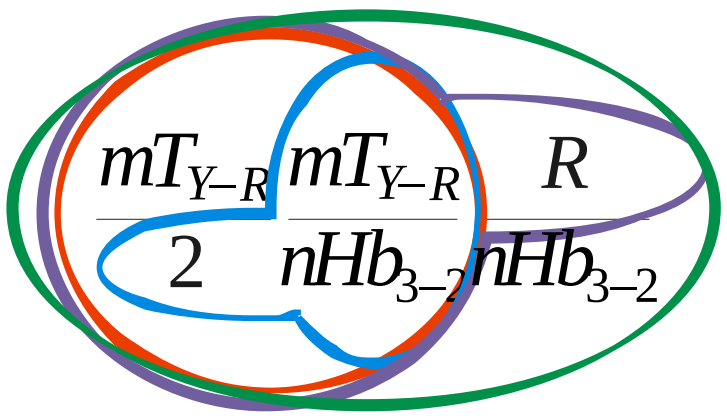

## 3 Concluding remarks

In our opinion, the set of rules needed to completely describe the genetic code degeneration allows to propose that a CD is an informative entity with a double function. The CD structure or composition informs about two closely related, but essentially different, phenomena. The rules here presented demonstrate that, apart from determining a specific AA, there is a “logic” that associates the CD composition with its discriminating or non-discriminating character. Together, these two aspects constitute the CD global information. Two different set of rules are needed to describe these different phenomena.

The rules here presented describe a CD as an informative entity whose global information depends on the spatially organized combinatory of the six properties conferred by its 3 constituent bases. In other words, as informative entity, every CD corresponds to a spatially organized combinations of four properties (**Y-R-2-3**) grouped as set of six elements. The spatial organiza-tion is emphasizes because the informative value of the properties Y, R, 2 and 3 is essentially linked to their specific position within the CD.

These rules establish that, depending on the kind of base (coherent or non-coherent) occupying the 2nd position, only two or three of six CD properties located at defined positions determine the discriminating or non-discriminating behavior. Besides, the number of properties, and their location within the CD, characteristically changes when an AA is specified by different kind of synonyms (4pts, 2pts, 1pts).

With regards to the discriminating or non-discriminating character of the CDs, it is interesting to remark the different informative role the base properties possess depending on their position within the CD structure. It is clear that the 2nd position possesses a relevant function: a coherent B2 *per se* determine the discriminating or non-discriminating character of CDs. In the case of CDs with non-coherent B2, however, additional information is required to specify the Discriminating vs non-Discriminating character. In considering a CD with a non-coherent B2, the other properties possess different roles depending on their location at B1 or at B3: **(a)** the property nHb of B1 is needed to determine the discriminating vs non-discriminating character and the property 2 allows the divergence of a 4pt into two 2pts. In this case, the property mT of B1 does not have informative value; **(b)** on the other hand, the property mT of B3 is needed to allow the divergence of a 2pt into 1pts. The property R in B3 allows the divergence of a 2pt into two 1pts while the property Y does not. Conversely, in this case the property nHb of B3 does not have informative value. On these bases it can be postulated that **a CD is an asymmetrically organized informative entity** since, clearly, the same informative elements possesses essentially different informative role when they are located in different positions, i.e. B1 or B3, with respect to a non-coherent B2.